\documentclass[preprint,sort&compress,final,11pt]{elsarticle}

\usepackage{amsmath}
\usepackage{amssymb}
\usepackage{graphicx}
\usepackage{latexsym}
\usepackage{epstopdf}
\usepackage{natbib}
\usepackage{color}
\usepackage{hyperref}

\newtheorem{theorem}{Theorem}[section]

\newtheorem{lemma}[theorem]{Lemma}

\numberwithin{equation}{section}

\def\Proof{\noindent{\bf Proof.}~}
\def\qed{\hfill$\square$\smallskip}

\def\dlim{\displaystyle\lim}

\journal{\empty}
\date{}

\begin{document}

\begin{frontmatter}

\title{Periodic solutions of sublinear impulsive differential equations }

\author[au1]{Yanmin Niu}

\ead[au1]{nym1216@163.com}

\author[au1]{Xiong Li\footnote{ Partially supported by the NSFC (11571041)and the Fundamental Research Funds for the Central Universities. Corresponding author.}}

\address[au1]{School of Mathematical Sciences, Beijing Normal University, Beijing 100875, P.R. China.}

\ead[au1]{xli@bnu.edu.cn}

\begin{abstract}
In this paper, we consider sublinear second order differential equations with impulsive effects. Basing on the Poincar\'{e}-Bohl fixed point theorem, we first will prove the existence of harmonic solutions. The existence of subharmonic solutions is also obtained by a new twist fixed point theorem recently established by Qian etc in 2015 (\cite{Qian15}).
\end{abstract}


\end{frontmatter}

\section{Introduction}
We are concerned in this paper with the existence of periodic solutions for the sublinear impulsive differential equation
\begin{equation}\label{maineq}
\left\{\begin{array}{ll}
 x''(t)+g(x)=p(t,x,x'),\quad t\neq t_{j};\\[0.2cm]
 \Delta x|_{t=t_{j}}=a x\left(t_{j}-\right),\\[0.2cm]
 \Delta x'|_{t=t_{j}}=a x'\left(t_{j}-\right),\quad\ \ \ \ \ j=\pm1,\pm2,\cdots,\\
 \end{array}
 \right.
 \end{equation}
where $0\leq t_{1}<\cdots<t_{k}<2\pi$, $a>0$ is a constant, $\Delta x|_{t=t_{j}}=x\left(t_{j}+\right)-x\left(t_{j}-\right)$ and $\Delta x'|_{t=t_{j}}=x'\left(t_{j}+\right)-x'\left(t_{j}-\right)$. In addition, we assume that the impulsive time is $2\pi$-periodic,
that is, $t_{j+k}=t_{j}+2\pi$ for $j=\pm1,\pm2,\cdots$, and $g$ is a continuous function with the sublinear growth condition
$$(g_{0}):\ \displaystyle\lim_{\mid x\mid\to+\infty}\frac{g(x)}{x}=0.$$

This problem comes from Duffing's equation
\begin{equation}\label{Duffing}
x''(t)+g(x)=p(t),
\end{equation}
and there is a wide literature dealing with this problem not only because its physical significance, but also the application of various mathematical techniques on it, such as Poincar\'{e}-Birkhoff twist theorem in \cite{Dingsuperquadratic}, \cite{Hartman77}, \cite{Jacobowitz76}, variational method in \cite{Bahri841}, \cite{Long89} and topological degree or index theories in \cite{Capietto92}, \cite{Capietto90}. Under different assumptions on the function $g$, for example being superlinear, subliner, superquadratic potential and so on, there are many interesting results on the existence and multiplicity of periodic solutions of \eqref{Duffing}, see \cite{Dingsublinear}, \cite{Pino}, \cite{Qian01}, \cite{Qian052}, \cite{Rebelo} and the references therein.

Recently, as impulsive equations widely arise in applied mathematics, they attract a lot of attentions and many authors study the general properties of them in \cite{Bainov93}, \cite{Lak}, along with the existence of periodic solutions of impulsive differential equations via fixed point theory in \cite{Nieto97}, \cite{Nieto02}, topological degree theory in \cite{Dong01}, \cite{Qian05}, and variational method in \cite{Nieto09}, \cite{Sun11}. However, different from the extensive study for second order differential equations without impulsive terms, there are only a few results on the existence and multiplicity of periodic solutions for impulsive second order differential equations.

In \cite{Qian15}, Qian etc considered the superlinear impulsive differential equation
\begin{equation}\label{qianeq}
\left\{
\begin{array}{ll}
 x''(t)+g(x)=p(t,x,x'),\quad \quad\quad\ \  t\neq t_{j};\\[0.2cm]
\Delta x|_{t=t_{j}}=I_{j} (x\left(t_{j}-\right),x'\left(t_{j}-\right)),\\[0.2cm]
\Delta x'|_{t=t_{j}}=J_{j} (x\left(t_{j}-\right),x'\left(t_{j}-\right)),\ \   j=\pm1,\pm2,\cdots,
 \end{array}
 \right.
\end{equation}
where $0\leq t_{1}<\cdots<t_{k}<2\pi$, $\Delta x|_{t=t_{j}}=x\left(t_{j}+\right)-x\left(t_{j}-\right)$, $\Delta x'|_{t=t_{j}}=x'\left(t_{j}+\right)-x'\left(t_{j}-\right)$,
and $I_{j}, J_{j}:\mathbb{R}\times\mathbb{R}\rightarrow\mathbb{R}$ are continuous maps for $j=\pm1,\pm2,\cdots$. In addition, they assume that the impulsive time is $2\pi$-periodic, that is, $t_{j+k}=t_{j}+2\pi$ for $j=\pm1,\pm2,\cdots$, and $g$ is a continuous function with the superlinear growth condition
$$(g_{0}'):\ \displaystyle\lim_{\mid x\mid\to\infty}\frac{g(x)}{x}=+\infty.$$ The authors proved via the Poincar\'{e}-Birkhoff twist theorem the existence of infinitely many periodic solutions of \eqref{qianeq} with $p=p(t)$ and also the existence
of periodic solutions for non-conservative case with degenerate impulsive terms by developing a new twist fixed point theorem.

In this article, we discuss \eqref{maineq} with the sublinear condition $(g_{0})$, which is different from the superlinear case and there is few papers studying on it up to now. As we all know, the existence of impulses, even the most simple impulsive functions, may cause complicated dynamic phenomenons and bring difficulties to study. Here, we start with linear impulsive functions and investigate the existence of periodic solutions completely. The rest is organized as follows. In Section 2, we obtain the existence of harmonic solutions of \eqref{maineq} by Poincar\'{e}-Bohl fixed point theorem. In Section 3, some properties of Poincar\'{e} map of \eqref{maineq}, including its derivative and rotation property, are considered and then the existence of subharmonic solutions is obtained according to a fixed point theorem coming from \cite{Qian15}.

\section{Harmonic solutions}
It is well known that Duffing's equation \eqref{Duffing} with the sublinear condition has at least one harmonic solution, see \cite{Mawhin72}. We will find that the existence of harmonic solutions also remains under the influence of special impulses, such as impulsive differential equation \eqref{maineq}.

\begin{theorem}\label{harmonicthe}
Suppose that in \eqref{maineq} the external force $p$ is a bounded and continuous function with $2\pi$-periodic  in the first variable, and $g$ is a locally Lipschitz continuous function satisfying the following conditions
$$(g_{0}):\ \displaystyle\lim_{\mid x\mid\to+\infty}\frac{g(x)}{x}=0$$ and
$$(g_{1}):\ \displaystyle\lim_{\mid x\mid\to+\infty}g(x)\mbox{sign}(x)=+\infty,\quad \displaystyle\lim_{x^{+}\to +\infty}\int_{x^{-}}^{x^{+}}\frac{ds}{\sqrt{G(x^{+})-G(s)}}=+\infty,$$
where $x^{-}<0<x^{+}$ are two zeros of the potential $G(x)=\int_{0}^{x}g(s)ds$. Then there is at least one $2\pi$-periodic solution of \eqref{maineq}.
\end{theorem}

Before giving the proof of Theorem \ref{harmonicthe}, we do some preparations. At first, some basic properties of impulsive differential equations will be given.

Consider the following initial value problem
\begin{equation}\label{generaleq}
\left\{
\begin{array}{ll}
 u'=f(t,u),\quad\ \ \ \ \ \ \ \  t\neq t_{j};\\[0.2cm]
 \Delta u|_{t=t_{j}}=I_{j}(u\left(t_{j}-\right)),\ j=\pm1,\pm2,\cdots,\\[0.2cm]
u(t_{0}+)=u_{0}, \\
 \end{array}
 \right.
\end{equation}
where $\Delta u|_{t=t_{j}}=u\left(t_{j}+\right)-u\left(t_{j}-\right),\ j=\pm1,\pm2,\cdots,$ and assume that
\begin{enumerate}
\item $f:\mathbb{R}\times\mathbb{R}^{n}\to\mathbb{R}^{n}$ is continuous in $(t_{j},t_{j+1}]\times\mathbb{R}^{n}$, locally Lipschitz in the second variable and the limits $\dlim_{t\to t_{j}+,v\to u}f(t,v),\ j=\pm1,\pm2,\cdots$ exist;
\item $I_{j}:\mathbb{R}^{n}\to\mathbb{R}^{n},\ j=\pm1,\pm2,\cdots,$ are continuous;
\item $f$ is $2\pi$-periodic in the first variable, $0\leq t_{1}<\cdots<t_{k}<2\pi$, $t_{j+k}=t_{j}+2\pi$ and $I_{j+k}=I_{j}$ for $j=\pm1,\pm2,\cdots$.
\end{enumerate}
Then the following lemma holds.

\begin{lemma}\label{generallem} (\cite{Bainov93}) Assume that the conditions above hold. Then for any $x_{0}\in\mathbb{R}^{n}$, there is a unique solution $u(t)=u(t;t_{0},u_{0})$ of \eqref{generaleq}. Moreover, $P_{t}:u_{0}\to u(t;t_{0},u_{0})$ is continuous in $u_{0}$ for $t\neq t_{j},\ j=\pm1,\pm2,\cdots$.
\end{lemma}

According to Lemma \ref{generallem}, we will further state some properties of the solutions of \eqref{generaleq} and $P_{t}$ in the next lemma, the proof can be found in \cite{Qian15}.

\begin{lemma}\label{generallem2}
Assume that the conditions of Lemma \ref{generallem} hold. Then all solutions of \eqref{generaleq} exist for $t\in\mathbb{R}$ provided that all solutions of $u'=f(t,u)$ exist for $t\in\mathbb{R}$. Moreover, if $\Phi_{j}: u\to u+I_{j}(u),\ j=1,2,\cdots,k$ are global homeomorphisms of $\mathbb{R}^{n}$, then $P_{t}$ is a homeomorphism for $t\neq t_{j},\ j=\pm1,\pm2,\cdots$. Furthermore, all solutions of \eqref{generaleq} have elastic property, that is, for any $b>0$, there is $r_{b}>0$ such that the inequality $|u_{0}|\geq r_{b}$ implies $|u(t;t_{0},u_{0})|\geq b$ for $t\in(t_{0},t_{0}+2\pi]$.
\end{lemma}

In what follows, we write \eqref{maineq} as an equivalent system of the form
\begin{equation}\label{maintwoeq}
\left\{
\begin{array}{ll}
 x'=y, \ \ y'=-g(x)+p(t,x,y),\ t\neq t_{j};\\[0.2cm]
 \Delta x|_{t=t_{j}}=a x\left(t_{j}-\right),\\[0.2cm]
 \Delta y|_{t=t_{j}}=a y\left(t_{j}-\right),\ \ \ \ \ \ \ \ \ \ j=\pm1,\pm2,\cdots.
 \end{array}
 \right.
\end{equation}
Let $x=x(t)=x(t;x_{0},y_{0}),\ y=y(t)=y(t;x_{0},y_{0})$ be the solution of \eqref{maintwoeq} satisfying the initial value condition $x(0)=x_{0}, y(0)=y_{0}$. By Lemma \ref{generallem2}, the solution $(x(t), y(t))$ exists for all $t\in\mathbb{R}$, and for any $b>0$, there is $r_{b}$ such that the inequality $r_{0}=\sqrt{x_{0}^{2}+y_{0}^{2}}\geq r_{b}$ implies $|r(t)|=\sqrt{x(t)^{2}+y(t)^{2}}\geq b$ for $t\in[0,2\pi]$. Denote by
$$P:(x_{0}, y_{0})\mapsto \left(x(2\pi;x_{0},y_{0}), y(2\pi;x_{0},y_{0})\right)$$
the Poincar\'{e} map of \eqref{maintwoeq}, then it is continuous in $(x_{0}, y_{0})$ and its fixed points are harmonic solutions of \eqref{maintwoeq} correspondingly.

By polar coordinates $$x=r\cos\theta, y=r\sin\theta,$$ \eqref{maintwoeq} has the following form
\begin{equation}\label{mainpolareq}
\left\{
\begin{array}{ll}
 \theta'=-\sin^{2}\theta-\left(g(r\cos\theta)-p(t,r\cos\theta,r\sin\theta)\right)\cos\theta/r,\\[0.2cm]
  r'=r\cos\theta\sin\theta-\left(g(r\cos\theta)-p(t,r\cos\theta,r\sin\theta)\right)\sin\theta,\quad t\neq t_{j};\\[0.2cm]
 \Delta \theta|_{t=t_{j}}=0,\\[0.2cm]
 \Delta r|_{t=t_{j}}=a r\left(t_{j}-\right), \ \ \ \ \ \ j=\pm1,\pm2,\cdots,\\
 \end{array}
 \right.
\end{equation}
from which we can find $\theta(t)$ is continuous in $t\in\mathbb{R}$ actually. Next in the phase plane, we give some general descriptions of the trajectory of \eqref{maintwoeq}.

\begin{lemma}\label{derivativelem} Assume that the conditions of Theorem \ref{harmonicthe} hold. If the initial value $r_{0}$ is large enough, then
$$\theta'(t;\theta_{0},r_{0})<0,\ t\in[0,2\pi],\ t\neq t_{j},\ j=1,2,\cdots,k$$ with $(\theta_{0},r_{0})=(\theta(0),r(0))$.
\end{lemma}
\Proof From the discussion above, we know that for any $c>0$, there is $d=d(c)>0$ such that for $t\in[0,2\pi]$, \begin{equation}\label{req}
r(t)=\displaystyle\sqrt{x(t)^{2}+y(t)^{2}}\geq c,
\end{equation}
if $r_{0}=\sqrt{x_{0}^{2}+y_{0}^{2}}\geq d$. Suppose that \eqref{req} holds and $c$ is a large number to be determined later.

In \eqref{mainpolareq}, we consider the first equation
\begin{equation}\label{polareq}
\theta'=-\sin^{2}\theta-\left(g(r\cos\theta)-p(t,r\cos\theta,r\sin\theta)\right)\cos\theta/r, t\neq t_{j},j=1,2,\cdots k.
\end{equation}
Due to the condition $(g_{1})$ and the boundedness of $p$, there exists a constant $N>0$ such that
\begin{equation}
\left\{
\begin{array}{ll}
g(x)-p(t,x,y)>\frac{1}{2}g(x),\ x\geq N;\\[0.2cm]
g(x)-p(t,x,y)<\frac{1}{2}g(x),\ x\leq -N,\\[0.2cm]
 \end{array}
 \right.
\end{equation}
and $x^{-1}g(x)>0$ for $|x|\geq N$. Then
\begin{equation}\label{negative1}
\theta'<-\sin^{2}\theta-\frac{1}{2}x^{-1}g(x)\cos^{2}\theta<0,\ |x|\geq N.
\end{equation}
On the other hand, when $|x|\leq N$, since $g$ and $p$ are bounded, there is a $K>0$ such that
$$-[g(x)-p(t,x,y)]\cos\theta\leq K, \ t\in[0,2\pi],$$
which implies that for $t\in[0,2\pi]$ and $|x|\leq N$,
\begin{equation}\label{negativesemi}
\theta'\leq K/r-\sin^{2}\theta.
\end{equation}
We choose $c>N$ large enough in \eqref{req} satisfying $\frac{\pi}{4}<\arccos\frac{N}{c}<\frac{\pi}{2}$. Denote $$\alpha=\arccos\frac{N}{c},$$ then for $r(t)>c$ and $|x|\leq N$, we have $|\cos\theta|\leq N/c=\cos\alpha$ and
$$\sin^{2}\theta\geq\sin^{2}\alpha\geq\frac{1}{2},$$
which implies
\begin{equation}\label{negative2}
\theta'\leq K/r-\frac{1}{2}<0,
\end{equation}
 if $c>2K$.

Combining \eqref{negative1} and \eqref{negative2} yields $\theta'(t;\theta_{0},r_{0})<0,\ t\in[0,2\pi],\ t\neq t_{j},\ j=1,2,\cdots,k$, for $r_{0}$ large enough. Thus we have finished the proof of the lemma.
\qed

\begin{lemma}\label{no2pilem}
Assume that the condition $(g_{1})$ holds. Then
\begin{equation}
-2\pi<\theta(2\pi;\theta_{0},r_{0})-\theta(0;\theta_{0},r_{0})<0
\end{equation}
with the initial value $r_{0}$ large enough.
\end{lemma}
\Proof Denote $\Gamma_{z_{0}}:z(t;z_{0})=(x(t;z_{0}),y(t;z_{0}))$ the trajectory of \eqref{maintwoeq} with the initial value condition $z_{0}=(x_{0},y_{0})$ and $(\theta(t;\theta_{0},r_{0}),r(t;\theta_{0},r_{0}))$ the polar coordinates of $z(t;z_{0})$.

As we know for $c_{0}>0$, there is $d_{0}>0$ such that the inequality $|z_{0}|\geq d_{0}$ implies $| z(t;z_{0})|\geq c_{0}$ for $t\in[0,2\pi]$. Let $d_{0}$ large and by Lemma \ref{derivativelem}, one has $\Gamma_{z_{0}}\in D$, where
$$
D=\{(x,y)\in\mathbb{R}^{2}:|(x,y)|=\sqrt{x^{2}+y^{2}}\geq c_{0}\},
$$
and
$$
\theta'(t;\theta_{0},r_{0})<0, \ t\in[0,2\pi],\ t\neq t_{j},\ j=1,2,\cdots k.
$$

In what follows, we divide $D$ into six subregion
$$
\begin{array}{ll}
&D_{1}=\{(x,y)\in D: \ |x|\leq c_{0},y>0\};\\[0.2cm]
&D_{2}=\{(x,y)\in D: \ x\geq c_{0},y\geq0\};\\[0.2cm]
&D_{3}=\{(x,y)\in D: \ x\geq c_{0},y\leq0\};\\[0.2cm]
&D_{4}=\{(x,y)\in D: \ |x|\leq c_{0},y<0\};\\[0.2cm]
&D_{5}=\{(x,y)\in D: \ x\leq -c_{0},y\leq0\};\\[0.2cm]
&D_{6}=\{(x,y)\in D: \ x\leq -c_{0},y\geq0\}.
\end{array}
$$
Since for $t\in[0,2\pi]$ and $t\neq t_{j}$, $\theta'(t)<0$ and for $t=t_{j}$, $\Delta\theta|_{t=t_{j}}=0,\ j=1,2,\cdots k$, $\Gamma_{z_{0}}$ moves by the following order:
$$D_{1}\rightarrow D_{2}\rightarrow D_{3}\rightarrow D_{4}\rightarrow D_{5}\rightarrow D_{6}\rightarrow D_{1}.$$

We first verify
$$
-2\pi<\theta(t_{j}-;\theta_{0},r_{0})-\theta(t_{j-1}+;\theta_{0},r_{0})<0, \ \ j=1,2,\cdots k,k+1,
$$
where $t_{0}=0$ and $t_{k+1}=2\pi$. If it is not true, $\Gamma_{z_{0}}$ moves at least one turn around the origin in $D$
when $t\in(t_{j-1},t_{j})$ for some $j=1,2,\cdots k,k+1$. Therefore, there exists at least an interval $[s_{1},s_{2}]\subset(t_{j-1},t_{j})$, during which $\Gamma_{z_{0}}$ goes through the region $D_{2}$. For the fixed $c_{0}$ and $d_{0}$ mentioned above, choose $d>d_{0}$ and there is a large constant $R_{0}$ such that the inequality $|z_{0}|>d>d_{0}$ implies $|z(t)|\geq R_{0},\ t\in[0,2\pi]$. When $t\in(s_{1},s_{2})$,
$$
x(t)>c_{0},\ \frac{dx(t)}{dt}=y(t)>0,
$$
and
$$
x_{1}=x(s_{1})=c_{0},\ y(s_{2})=0,\ x_{2}=x(s_{2})>R_{0}.
$$
Since $p$ is bounded, denote
$$
m=\inf_{(t,x,y)\in[0,2\pi]\times\mathbb{R}\times\mathbb{R}}|p(t,x,y)|,
$$
and then by \eqref{maintwoeq},
$$y(t)\displaystyle\frac{dy(t)}{dt}=-g(x(t))\displaystyle\frac{dx(t)}{dt}+p(t,x(t),y(t))\displaystyle\frac{dx(t)}{dt}\geq-(g(x(t))-m)\displaystyle\frac{dx(t)}{dt}.
$$
Integrating the inequality above on the interval $[t,s_{2}]\subset[s_{1},s_{2}]$, we have
$$\begin{array}{lll}
\frac{1}{2}y^{2}(t)&\leq&G(x_{2})-G(x(t))-m(x_{2}-x(t))\\
&=&[G(x_{2})-G(x(t))]\left(1-\displaystyle\frac{m}{g(\xi)}\right),\\
\end{array}
$$
where $\xi\in(x_{1},x_{2})$ by the mean value theorem. From the condition $(g_{1})$, one has
$$0<\displaystyle\left(1-\displaystyle\frac{m}{g(\xi)}\right)<1$$ for $R_{0}$ large enough. Therefore
$$\frac{dx(t)}{dt}=y(t)\leq\sqrt{2(G(x_{2})-G(x(t)))},\quad s_{1}<t<s_{2},$$
from which we deduce
$$s_{2}-s_{1}\geq\displaystyle\int_{c_{0}}^{x_{2}}\frac{dx}{\sqrt{2}\cdot\sqrt{(G(x_{2})-G(x(t)))}}.$$
Notice that $c_{0}$ is a fixed constant and $x_{2}>R_{0}$. By the condition $(g_{1})$, for $R_{0}$ sufficient large, $s_{2}-s_{1}$ can be large enough, which implies $t_{j}-t_{j-1}$ cannot be bounded for $j=1,2,\cdots k,k+1$. That is a contradiction.
Then $$-2\pi<\theta(t_{j}-;\theta_{0},r_{0})-\theta(t_{j-1}+;\theta_{0},r_{0})<0, j=1,2,\cdots k,k+1,$$ where $t_{0}=0$ and $t_{k+1}=2\pi$.

Finally, from the proof above, together with the condition $\Delta\theta|_{t=t_{j}}=0,\ j=1,2,\cdots k$, $\Gamma_{z_{0}}$ cannot go through the region $D_{2}$ for $t\in[0,2\pi]$. Therefore $$-2\pi<\theta(2\pi;\theta_{0},r_{0})-\theta(0;\theta_{0},r_{0})<0$$ with the initial value $r_{0}$ large enough.
\qed

We are now in a position to prove Theorem \ref{harmonicthe}.\\

\noindent {\bf Proof of Theorem \ref{harmonicthe}} By Lemma \ref{no2pilem}, for $|z_{0}|=d$ large enough, the argument of $\Gamma_{z_{0}}$ satisfies $$-2\pi<\theta(2\pi;z_{0})-\theta(0;z_{0})<0.$$ Then on $D_{d}\triangleq\{z\in\mathbb{R}^{2}:| z|\leq d\}$, the Poincar\'{e} map $P$ of \eqref{maintwoeq} has at least one fixed point $\zeta_{0}$ by Poincar\'{e}-Bohl fixed point theorem. Correspondingly, $z(t;\zeta_{0})$ is the $2\pi$-periodic solution of \eqref{maintwoeq}.\qed

\section{Subharmonic solutions}

In this section, we discuss the existence of subharmonic solutions of \eqref{maineq}, and we always assume that the conditions $(g_{0})$, $(g_{1})$ hold and
$$
(g_{2}): g'(x)>0, \ x\in \mathbb{R}.
$$
We now state the existence of subharmonic solutions of \eqref{maineq}.
\begin{theorem}\label{subthe}
Assume that the conditions above hold. Then there is a positive integer sequence $$0<n_{1}<n_{2}<\cdots<n_{k}<\cdots (\rightarrow \infty)$$
such that \eqref{maineq} has at least one $2n_{k}\pi$-periodic solution for each $n_{k},\ k=1,2,\cdots$.
\end{theorem}

By Theorem \ref{harmonicthe}, Eq.\eqref{maintwoeq} has at least one harmonic solution: $x=\overline{x}(t),\ y=\overline{y}(t)$. Let
$$x=u+\overline{x}(t),\  y=v+\overline{y}(t),$$
then \eqref{maintwoeq} can be transformed into
\begin{equation}\label{uvtwoeq}
\left\{
\begin{array}{ll}
 u'=v, \ \ v'=-H(t,u)u,\quad t\neq t_{j};\\[0.2cm]
 \Delta u|_{t=t_{j}}=a u\left(t_{j}-\right),\\[0.2cm]
 \Delta v|_{t=t_{j}}=a v\left(t_{j}-\right), \quad\ \ \ \ \ j=\pm1,\pm2,\cdots,
 \end{array}
 \right.
\end{equation}
where $$H(t,u)=\displaystyle\int_{0}^{1}g'\left(\overline{x}(t)+su\right)ds$$ is $2\pi$-periodic in the first variable. By $(g_{2})$, $H(t,u)>0$.

It is easy to see that $(u,v)=(0,0)$ is a trivial solution of \eqref{uvtwoeq}. Let $(u_{0},v_{0})=(u(0),v(0))\neq0$, then the solution $(u(t),v(t))$ of \eqref{uvtwoeq} with the initial value $(u_{0},v_{0})$ cannot intersect with the origin. Denote the solution by polar coordinates $$u=\rho\cos\varphi,\ v=\rho\sin\varphi.$$
If $\rho_{0}=\sqrt{u_{0}^{2}+v_{0}^{2}}\neq 0$, then $\rho=\rho(t)=\rho(t;\rho_{0},\varphi_{0})>0$, where $\varphi_{0}$ is the argument of $(u_{0},v_{0})$ with $\varphi_{0}\in[0,2\pi)$. From \eqref{uvtwoeq}, it follows that
$$\Delta \varphi|_{t=t_{j}}=0,\  j=\pm1,\pm2,\cdots,$$
which implies that $\varphi=\varphi(t)=\varphi(t;\rho_{0},\varphi_{0})$ is continuous in $(t,\rho_{0},\varphi_{0})\in \mathbb{R}\times \mathbb{R}^{+}\times[0,2\pi)$.

Under polar coordinates,
\begin{equation}\label{uvpolareq}
\varphi'=-(H(t,\rho\cos\varphi)\cos^{2}\varphi+\sin^{2}\varphi), \ t\neq t_{j},\ j=\pm1,\pm2,\cdots.
\end{equation}
Then for $\rho_{0}>0$,
\begin{equation}\label{uvderivativeeq}
\varphi'(t)<0,\ t\in \mathbb{R},\ t\neq t_{j},\ j=\pm1,\pm2,\cdots.
\end{equation}

Next, we prove the following lemma.
\begin{lemma}\label{jiaodu}
The equality
\begin{equation}\label{inftyeq}
\displaystyle\lim_{t\rightarrow +\infty}\varphi(t)=-\infty
\end{equation}
holds for $\rho_{0}>0$.
\end{lemma}
\Proof We first assume that $\rho=\rho(t)$ is bounded as $t\to+\infty$. Since $g'\left(\overline{x}+su\right)$ is bounded on $s\in[0,1]$, $H$ has a minimum $H_{0}$ such that $$H(t,u)\geq H_{0}>0, \ t\geq 0.$$ By \eqref{uvpolareq}, we have the estimate $$\varphi'(t)\leq-(H_{0}\cos^{2}\varphi+\sin^{2}\varphi)\leq-a,$$ where $a\triangleq\displaystyle\inf_{\varphi\in \mathbb{R}}(H_{0}\cos^{2}\varphi+\sin^{2}\varphi)>0$. That implies \eqref{inftyeq}.

Secondly, we assume that $\rho=\rho(t)$ is unbounded as $t\to+\infty$. It follows from \eqref{uvderivativeeq} that $$\displaystyle\lim_{t\rightarrow +\infty}\varphi(t)=\varphi_{1}\geq-\infty.$$ Then we just need to prove $\varphi_{1}=-\infty$. Otherwise, we assume that $\varphi_{0}>\varphi_{1}>-\infty$. Firstly, by \eqref{uvpolareq}, we have $\varphi_{1}=k\pi$ ($k$ is some integer). Indeed, if $\varphi_{1}\neq k\pi$, then $\sin^{2}\varphi_{1}\neq0$. By \eqref{uvpolareq}, for large enough $t$,  we have
$$\varphi'(t)<-\sin^{2}\varphi(t)<-\frac{1}{2}\sin^{2}\varphi_{1},$$
which means $\varphi(t)\rightarrow-\infty$ as $t\rightarrow+\infty$, and contradicts with the assumption.

Without loss of generality, let $k=0$ and
$$\displaystyle\lim_{t\rightarrow +\infty}\varphi(t)=0.$$
Then for $(u_{0},v_{0})$ belonging to the first quadrant, $$\frac{du}{dt}=v>0,\ \Delta u|_{t=t_{j}}=a u\left(t_{j}-\right)>0,$$
from which we know that
$$
\displaystyle\lim_{t\to +\infty}u(t)=+\infty
$$
and
$$
\displaystyle\lim_{t\to +\infty}v(t)=0+.
$$
We notice that $(\overline{x}(t),\overline{y}(t))$ is the harmonic solution of \eqref{maintwoeq}, therefore
$$
\displaystyle\lim_{t\to +\infty}x(t)=\displaystyle\lim_{t\to +\infty}(u(t)+\overline{x}(t))=+\infty.
$$
By \eqref{maintwoeq} and the condition $(g_{1})$, as well as the boundedness of $p$, one has for $t$ large enough, there is a constant $M>0$ such that
$$
y'(t)=(-g(x(t))+p(t,x,y))<-M,
$$
which implies
\begin{equation}\label{ynegfintyeq}
\displaystyle\lim_{t\to +\infty}y(t)=-\infty.
\end{equation}
On the other hand, since $\displaystyle\lim_{t\to +\infty}v(t)=0+$, it follows that
$y(t)=v(t)+\overline{y}(t)$ is bounded as $t\rightarrow+\infty$, which contradicts with \eqref{ynegfintyeq}. Thus $\varphi_{1}=-\infty$ and the proof is finished. \qed

Generally speaking, if the external force $p$ in \eqref{Duffing} depends on the derivative $x'$, then equation \eqref{Duffing} is not conservative, and hence the Poincar\'{e} map of \eqref{Duffing} is not area-preserving and the method via the Poincar\'{e}-Birkhoff twist theorem is invalid. To overcome this difficulty, Qian etc \cite{Qian15} established a new fixed point theorem, which is a partial extension of the Poincar\'{e}-Birkhoff twist theorem. As for the non-conservative equation \eqref{maineq}, we will find that this theorem can be used to obtain the existence of subharmonic solutions. The fixed point theorem is stated as follows.

\begin{lemma}\label{fixlem}(\cite{Qian15})
Let $\Gamma_{-}$ and $\Gamma_{+}$ be two convex closed curves surrounding the origin, $int(\Gamma_{-})$ and $int(\Gamma_{+})$ be the interior domain of $\Gamma_{-}$ and $\Gamma_{+}$, respectively. Denote by $\mathcal{A}$ the annulus bounded by $\Gamma_{-}$ and $\Gamma_{+}$. Consider a continuous map $F:\overline{int(\Gamma_{+})}\rightarrow\mathbb{R}^{2}$. Let
$$E=\{z\in\mathcal{A}:|F(z)|\leq |z|\},$$
$$J=\{z\in\mathcal{A}:F(z)\in\mathbb{R}^{2}\backslash U(O),\ <Lz,F(z)>=0\}, $$
where $U(O)$ is a small neighborhood of the origin O and $L$ is a real orthogonal matrix with $det (L)=1$.

If for any curve $\gamma$ connecting $\Gamma_{-}$ and $\Gamma_{+}$, one has $\gamma\cap(J\cup E)\neq\emptyset$, then $F$ has at least one fixed point in $\overline{int(\Gamma_{+})}$.
\end{lemma}

In the $u,v$ plane, let $$B_{r}=\{(u,v)\in\mathbb{R}^{2}:u^{2}+v^{2}<r^{2}\},$$
$$C_{r}=\{(u,v)\in\mathbb{R}^{2}: u^{2}+v^{2}=r^{2}\},$$
and $\Lambda_{\omega_{0}}:w(t;w_{0})=(u(t),v(t))$ be the trajectory of \eqref{uvtwoeq} with the initial value $w_{0}=(u_{0},v_{0})=(u(0),v(0))$. By Lemma \ref{generallem2}, for
$b>0$, there is $c>0$ such that $|w_{0}|\geq c$ implies $|w(t)|\geq b,\ t\in[0,2\pi]$.
In order to give a description of the motion of $\Lambda_{\omega_{0}}$, we consider the region
$$\Omega:|\varphi|<\frac{\pi}{4},\ \ \rho\geq b, $$
and estimate the time of $\Lambda_{\omega_{0}}$ going through $\Omega$ with $w_{0}\in C_{c}$.

In $\Omega$, one has
$$u=\rho\cos\varphi\geq b \cos\frac{\pi}{4}=\frac{\sqrt{2}}{2}b.$$
According to the condition $(g_{0})$, for any $\delta>0$, there exists some $b$ large enough such that
\begin{equation}\label{smallheq}
H(t,u)=\frac{g(\overline{x}+u)-g(\overline{x})}{u}<\delta^{2},\ \ \mbox{if}\  u\geq\frac{\sqrt{2}}{2}b.
\end{equation}
Then in $\Omega$ and for $t\neq t_{j},\ j=1,2,\cdots,k$, by \eqref{uvpolareq} we have
$$\begin{array}{lll}
0>\varphi'&=&-\left(H(t,\rho\cos\varphi)\cos^{2}\varphi+\sin^{2}\varphi\right)\\[0.2cm]
&>&-\delta^{2}\cos^{2}\varphi-\sin^{2}\varphi.\\
\end{array}
$$
As for $t=t_{j},\ j=1,2,\cdots,k$, one has $\Delta \varphi|_{t=t_{j}}=0$ which implies that $\varphi(t)$ is continuous in $t\in \mathbb{R}$. Then the time of $\Lambda_{\omega_{0}}$ going through $\Omega$ outside $B_{b}$ is
$$\tau>\displaystyle\int_{\frac{\pi}{4}}^{-\frac{\pi}{4}}\frac{-d\varphi}{\delta^{2}\cos^{2}\varphi+\sin^{2}\varphi}=\frac{2}{\delta}\arctan\frac{1}{\delta}.$$
Because $$\displaystyle\lim_{\delta\to 0+}\frac{2}{\delta}\arctan\frac{1}{\delta}=+\infty,$$ then for any $n\in\mathbb{Z}^{+}$, we have
\begin{equation}\label{timeinftyeq}
\tau>2n\pi,
\end{equation}
if $\delta$ is sufficiently small (or the initial value $c$ mentioned above is sufficiently large).

\begin{lemma}\label{sublem}
Let the Poincar\'{e} map of \eqref{uvtwoeq} be $$P_{0}:(u_{0},v_{0})\rightarrow(u(2\pi;u_{0},v_{0}),v(2\pi;u_{0},v_{0})).$$ Then for any $n\in\mathbb{Z}^{+}$, there are $c_{n}>0$ and $n^{\ast}>n$ such that
$P_{0}^{k}\ (k\leq n)$ has no fixed point while $P_{0}^{n^{\ast}}$ has at least one fixed point outside $B_{c_{n}}$.
\end{lemma}
\Proof For any $n\in\mathbb{Z}^{+}$, we choose $c_{n}$ large enough so that \eqref{timeinftyeq} holds, which implies that $P_{0}^{k}\ (k\leq n)$ has no fixed point outside $B_{c_{n}}$.

Fix the $c_{n}$ above. By \eqref{inftyeq}, for $(u_{0},v_{0})\in C_{c_{n}}$, there is a positive integer $n^{\ast}>n$ such that
$$\varphi(2n^{\ast}\pi)-\varphi(0)<-2\pi+\frac{\pi}{2}.$$

On the other hand, \eqref{timeinftyeq} guarantees that there exists $e_{n} (e_{n}>c_{n})$ such that
$$-2\pi+\frac{\pi}{2}<\varphi(2n^{\ast}\pi)-\varphi(0)<0, \ (u_{0},v_{0})\in C_{e_{n}}.$$

Denote by $\mathcal{A}=\{(u,v)\in\mathbb{R}^{2}:\ c_{n}^{2}\leq u^{2}+v^{2}\leq e_{n}^{2}\}$ the annulus bounded by $C_{c_{n}}$ and $C_{e_{n}}$. Let $\beta$ be a curve connecting $C_{c_{n}}$ and $C_{e_{n}}$, then there is $\omega^{\ast}\in\beta$ such that $$\varphi(2n^{\ast}\pi;\varphi^{\ast},\rho^{\ast})-\varphi(0;\varphi^{\ast},\rho^{\ast})=-2\pi+\frac{\pi}{2},$$
where $(\varphi^{\ast},\rho^{\ast})$ are the polar coordinates of $\omega^{\ast}$. Therefore, choosing $L=id$ in Lemma \ref{fixlem}, we have $$<L\omega^{\ast},P_{0}^{n^{\ast}}(\omega^{\ast})>=0,$$
which implies $\beta\cap J\neq\emptyset$. The continuous map $P_{0}^{n^{\ast}}$ meets all assumptions of Lemma \ref{fixlem}. Thus $P_{0}^{n^{\ast}}$ has at least one fixed point $\omega_{\ast}$ in $\mathcal{A}$.\qed

Actually, Theorem \ref{subthe} is the corollary of Lemma \ref{sublem}, since the fixed point of $P_{0}^{n^{\ast}}$ is also the $2n^{\ast}\pi$-periodic solutions of \eqref{uvtwoeq}. Then we have finished the proof of Theorem \ref{subthe}.

\section*{References}
\bibliographystyle{elsarticle-num}

\end{document}